\documentclass[12pt]{amsart}
\newtheorem{thm}{Theorem}
\newtheorem{lem}{Lemma}
\newtheorem*{clm}{Claim}
\newtheorem*{scl}{Subclaim}
\newtheorem*{cor}{Corollary}
\newcommand{\al}{\alpha}
\newcommand{\be}{\beta}
\newcommand{\ga}{\gamma}

\newcommand{\ka}{\kappa}
\newcommand{\la}{\lambda}
\newcommand{\om}{\omega}
\newcommand{\si}{\sigma}

\newcommand{\cf}{\operatorname{cf}}
\newcommand{\cl}{\operatorname{cl}}
\newcommand{\cov}{\operatorname{cov}}
\newcommand{\su}{\subset}
\newcommand{\lan}{\langle}
\newcommand{\ran}{\rangle}
\newcommand{\p}{{\mathcal P}}
\newcommand{\pkl}{\p_\ka\la}
\newcommand{\lm}{[\la]^\mu}
\newcommand{\fo}{\forall}
\newcommand{\ex}{\exists}
\newcommand{\bi}{\bigcup}

\begin{document}
\title[\null]{Nonreflecting stationary sets in $\pkl$}
\author[\null]{Saharon Shelah and Masahiro Shioya}
\address{Institute of Mathematics,
Hebrew University,
Jerusalem, 91904 Israel.}
\email{shelah@math.huji.ac.il}
\address{Institute of Mathematics,
University of Tsukuba,
Tsukuba, 305-8571 Japan.}
\email{shioya@math.tsukuba.ac.jp}
\thanks{The first author was 
supported by the Israel Science
Foundation founded by
the Israel Academy of Sciences and Humanities.
Publication 764.
The second author was partially
supported by Grant-in-Aid for Scientific Research
(No.12640098),  
Ministry of Education, Science, Sports
and Culture of Japan.}
\subjclass{03E05}
\keywords{}
\begin{abstract}
Let $\ka$ be a regular uncountable cardinal
and $\la\ge\ka^+$.
The principle of stationary reflection for
$\pkl$ has been successful in settling
problems of infinite combinatorics
in the case $\ka=\om_1$.
For a greater $\ka$ 
the principle is known to fail
at some $\la$.
This note shows
that it fails at every $\la$ if
$\ka$ is the successor of a regular
uncountable cardinal
or $\ka$ is countably closed.
\end{abstract}
\date{}
\maketitle

\section{Introduction}
In \cite{fms} 
Foreman, Magidor and Shelah introduced
the following principle 
for $\la\ge\om_2$:
If $S$ is a stationary subset of 
$\p_{\om_1}\la$, then
$S\cap\p_{\om_1}A$ is stationary in 
$\p_{\om_1}A$
for some $\om_1\su A\su\la$
of size $\om_1$.
Let us call the principle stationary 
reflection for $\p_{\om_1}\la$.
It follows from Martin's Maximum
(see \cite{fms}) and holds
in the L\'evy model where 
$\om_2$ was supercompact 
in the ground model (see \cite{be}).
See \cite{cfm, s4, sh0, w}
for recent applications
of reflection principles for 
stationary sets in $\p_{\om_1}\la$.

What if $\om_1$ is replaced 
by a higher regular cardinal?
Feng and Magidor~\cite{fem} proved that
the corresponding statement for $\p_{\om_2}\la$
is false at some
large enough $\la$.
Their argument (see also \cite{be})
showed in effect that
stationary reflection for $\pkl$ 
at some large enough $\la$ implies 
the presaturation of
the club filter on $\ka$
for a successor cardinal $\ka$, which
is known to be false if in addition
$\ka\ge\om_2$ by \cite{s0}.

Foreman and Magidor~\cite{fom} 
extended the Feng--Magidor result
for every regular cardinal $\ka\ge\om_2$,
although they proved only the case
$\ka=\om_2$.
We present below what was  
proved in effect
and in \S4 its proof of our own:

\begin{thm} 
Stationary reflection for $\pkl$ fails
at every $\la\ge2^{\ka^+}$ if $\ka\ge\om_2$
is regular.
\end{thm}

See \cite{s2}
for a further example of nonreflection,
which is based on pcf theory~\cite{s1}.
This note
addresses the problem whether
stationary reflection for $\pkl$
fails {\it everywhere}, i.e. 
at every $\la\ge\ka^+$.
Specifically we prove

\begin{thm} 
Stationary reflection for $\pkl$ fails
everywhere if 
$\nu<\ka$ are both regular uncountable
and
$\cf(\nu,\ga)<\ka$ for $\nu<\ga<\ka$.
\end{thm}

Here
$\cf(\nu,\ga)$ is the smallest size
of unbounded subsets of $\p_\nu\ga$.
The last condition in Theorem~2 holds
if $\ka=\nu^+$ or if $\nu=\om_1$ and
$\ga^\om<\ka$ for $\ga<\ka$.
In \S3 we prove Theorem~2
in much greater generality.

\section{Preliminaries}

For background material
we refer the reader to \cite{k}.
Throughout the paper, 
$\ka$ and $\nu$ stand for
a regular cardinal $\ge\om_1$
and
$\mu<\la$ a cardinal $\ge\ka$.
We write $S^\nu_\ka$ for 
$\{\ga<\ka:\cf\ga=\nu\}$.
Let $A$ be a set of ordinals.
The set of limit points of
$A$ is denoted $\lim A$.
It is easy to see 
$|\lim A|\le|A|$.
$A$ is called $\si$-closed if
$\ga\in A$ for $\ga\in\lim A$ of cofinality
$\om$.
Let 
$f:[\la]^{<\om}\to\pkl$.
We write
$C(f)$ for 
$\{x\in\pkl:\bi f``[x]^{<\om}\su x\}$.
For $x\in\pkl$ 
the smallest superset of $x$
in $C(f)$ is denoted $\cl_fx$.

Stationary reflection for $\pkl$ states that
if $S$ is a stationary subset of $\pkl$, then
$S\cap\p_\ka A$ is stationary in $\p_\ka A$
for some $\ka\su A\su\la$ of size $\ka$.
It is easily seen that
stationary reflection for $\pkl$
implies one for $\p_\ka\mu$.
Hence stationary reflection for
$\pkl$ fails everywhere 
iff it fails at $\la=\ka^+$.

Let $S$ be a stationary subset of $\pkl$.
$S$ is called nonreflecting
if 
it witnesses the failure of
stationary reflection, i.e.
$S\cap\p_\ka A$ is nonstationary in $\p_\ka A$
for $\ka\su A\su\la$ of size $\ka$.
More generally
$S$ is called $\mu$-nonreflecting if
$S\cap\p_\ka A$ is nonstationary in $\p_\ka A$
for $\mu\su A\su\la$
of size $\mu$.

We write $[\la]^\mu$ for 
$\{x\su\la:|x|=\mu\}$.
A filter $F$ on $[\la]^\mu$ 
is called fine if it is $\mu^+$-complete
and
$\{x\in\lm:\al\in x\}\in F$
for $\al<\la$.
The specific example relevant to us was
introduced in \cite{s}:

\begin{lem}
A fine filter on $\lm$ is generated by
the sets of the form
$\{\bi_{n<\om}A_n:\{A_n:n<\om\}\su\lm\land
\fo n<\om(\varphi(\lan A_k:k<n\ran)\su A_n)\}$,
where $\varphi:(\lm)^{<\om}\to\lm$.
\end{lem}

We need an analogue
\cite{m} of Ulam's theorem in our context:

\begin{lem} 
$\lm$ splits into $\la$ disjoint 
$F$-positive sets 
if $F$ is a fine filter on
$\lm$. 
\end{lem}

\begin{proof}
It suffices to split
$X $ $F$-positive into
$\nu$ disjoint $F$-positive sets for
$\mu<\nu\le\la$ regular.
Fix a bijection
$\pi_x:\mu\to x$ for $x\in X$.
Set
$X_{\ga\xi}=\{x\in X:\pi_x(\xi)=\ga\}$
for
$\ga<\nu$ and $\xi<\mu$.
Then 
$\bi_{\xi<\mu}X_{\ga\xi}=
\{x\in X:\ga\in x\}$ is  $F$-positive
for $\ga<\nu$.
Hence for $\ga<\nu$
we have $\xi<\mu$ such that
$X_{\ga\xi}$ is $F$-positive,
since $F$ is $\mu^+$-complete.
Thus we have $F$-positive sets
$\{X_{\ga\xi}:\ga\in A\}\su\p X$ 
for some
$A\in[\nu]^\nu$
and $\xi<\mu$, which are mutually disjoint,
as desired.
\end{proof}

\section{Main Theorem}
This section is devoted to the main result
of this paper.
Like the proof \cite{s0} of 
a diamond principle for some $\p_{\om_1}\la$
(see also \cite{sh}),
our argument originates from 
nonstructure theory \cite{s3}.

Throughout the section, let
$\nu<\ka$ be regular cardinals $\ge\om_1$
and $\mu<\la$
cardinals $\ge\ka$.
Recall from \cite{s1} 
$\cov(\la,\mu^+,\mu^+,\nu)=\la$ iff
$\{\bi_{\al\in a}E_\al:a\in\p_\nu\la\}$ 
is unbounded in $\lm$ 
for some 
$\{E_\al:\al<\la\}\su[\la]^\mu$.
It is easy to see
$\cov(\mu^+,\mu^+,\mu^+,\nu)=\mu^+$.

For the moment assume further
$\cf(\nu,\ga)<\ka$ for $\nu<\ga<\ka$.
Inductively we have
$\{c_\xi:\xi<\ka\}\su\p_\nu\ka$
and
$g:\ka\to\ka$ so that
$\{c_\xi:\xi<g(\ga)\}$ is unbounded in
$\p_\nu\ga$ for $\nu\le\ga<\ka$. 
Then 
$T=\{\ga\in S^\nu_\ka:g``\ga\su\ga\}$
is stationary in $\ka$ and
$\{c_\xi:\xi<\ga\}$
is unbounded in $\p_\nu\ga$
for $\ga\in T$.
Hence Theorem~2 
follows from the case
$\la=\mu^+=\ka^+$ of

\begin{thm} 
Assume
$\cov(\la,\mu^+,\mu^+,\nu)=\la$, 
$\{c_\xi:\xi<\mu\}\su\p_\nu\mu$,
$T$ is a stationary subset of $\p_\ka\mu$
of size $\mu$ and
$\{c_\xi:\xi\in z\}$
is unbounded in $\p_\nu z$
for $z\in T$. 
Then $\pkl$ has
a $\mu$-nonreflecting stationary subset.
\end{thm}

\begin{proof}
Let
$\{E_\al:\al<\la\}\su[\la]^\mu$
witness
$\cov(\la,\mu^+,\mu^+,\nu)=\la$.
Define $e:\la\times\mu\to\la$ so that
$E_\al=e``\{\al\}\times\mu$.
Hence for $A\in\p_{\mu^+}\la$
we have 
$a\in\p_\nu\la$ with
$A\su e``a\times\mu$.
Let $F$ be the filter on $\lm$ 
as defined in Lemma~1.
Lemma~2 allows us to split 
$[\la]^\mu$ into $\mu$ disjoint
$F$-positive sets 
$\{X_z:z\in T\}$.

Set
$S=\{x\in\pkl:e``x\times(x\cap\mu)\su x\land
x\cap\mu\in T\land
\ex b\in\p_\nu x
(x\su e``b\times\mu=
e``x\times\mu\in X_{x\cap\mu})\}$.

\begin{clm}
$S$ is stationary in $\pkl$. 
\end{clm}

\begin{proof}
Fix $f:[\la]^{<\om}\to\pkl$.
We may assume
$e``x\times(x\cap\mu)\su x$
for $x\in C(f)$.
For
$z\in T$
consider the following game ${\mathcal G}(z)$
of length $\om$
between two players {\it I} and {\it II}:

At round $n$ {\it I} plays 
$\mu\su A_n\su\la$ of size $\mu$.
Then
{\it II} plays a triple of 
$b_n\in\p_\nu\la$, a bijection 
$\pi_n:\mu\to e``b_n\times\mu$ 
and $x_n\in C(f)$
such that
$b_n\su x_n=\pi_n``(x_n\cap\mu)$.
We further require
$A_n\su e``b_n\times\mu\su 
e``x_n\times\mu\su A_{n+1}$ and
$x_n\su x_{n+1}$.
Finally we let
{\it II} win 
iff $x_n\cap\mu=z$ for $n<\om$.

Set
$T'=\{z\in T:$ {\it II}
has no winning strategy in 
${\mathcal G}(z)\}$.

\begin{scl}
$T'$ is nonstationary in $\p_\ka\mu$.
\end{scl}

\begin{proof}
Suppose otherwise.
For $z\in T'$ we have a winning strategy 
$\tau_z$ for {\it I} in ${\mathcal G}(z)$, since
the game is closed for {\it II},
hence determined.
By induction on $n<\om$ build 
$b_n$, $\pi_n$ and $\{x^z_n:z\in T'\}$
so that
$\lan(b_n,\pi_n,x^z_n):n<\om\ran$
is a play of {\it II} in ${\mathcal G}(z)$
against $\tau_z$ as follows:

Since $|T'|\le|T|=\mu$, 
we have in $\om$ steps
$b_n\in\p_\nu\la$ such that
$\bi_{z\in T'}\tau_z(\lan(b_k,\pi_k,x^z_k):k<n\ran)
\su e``b_n\times\mu$,
$b_n\su e``b_n\times\mu$
and $ e``b_n\times\mu$ is closed under $f$.
Next fix a bijection 
$\pi_n:\mu\to e``b_n\times\mu$.
Note that 
$x^z_{n-1}=\pi_{n-1}``(x^z_{n-1}\cap\mu)\su
e``b_{n-1}\times\mu\su 
\tau_z(\lan(b_k,\pi_k,x^z_k):k<n\ran)
\su e``b_n\times\mu$
for $z\in T'$.
Hence we have
$x^z_{n-1}\cup b_n\su x^z_n
\su e``b_n\times\mu$
such that
$\pi_n``(x^z_n\cap\mu)=x^z_n\in C(f)$,
since
$b_n\su e``b_n\times\mu$
and
$e``b_n\times\mu$ is closed under $f$.
If possible, we further require
$x^z_n\cap\mu=z$, in which case we have
$x^z_n=\pi_n``z$.

Set $b=\bi_{n<\om}b_n\in\p_\nu\la$
and $E=e``b\times\mu\in[\la]^\mu$.
Then $b\su E$ by
$b_n\su e``b_n\times\mu$.
Since
$e``b_n\times\mu\su e``b_{n+1}\times\mu$
are closed under $f$, so is $E$.
Also
$\mu\su\bi_{z\in T'}\tau_z(\emptyset)\su
e``b_0\times\mu\su E$.
Since $T'$ is stationary in $\p_\ka\mu$,
we have
$b\su x\su E$
such that
$x\in C(f)$,
$\pi_n``(x\cap\mu)=x\cap e``b_n\times\mu$ 
for $n<\om$ and
$x\cap\mu\in T'$.

Set $z=x\cap\mu$.
Since
$\mu\su e``b_0\times\mu\su e``b_n\times\mu$,
it is easily seen that
$x\cap e``b_n\times\mu=\pi_n``z$
meets the requirements for
$x^z_n$.
Hence
$x^z_n=x\cap e``b_n\times\mu$ and
$x^z_n\cap\mu=x\cap\mu=z$ for $n<\om$.
Thus {\it II} wins 
against $\tau_z$ with the play
$\lan(b_n,\pi_n,x^z_n):n<\om\ran$, 
which contradicts
that $\tau_z$ is a winning strategy for {\it I}
in ${\mathcal G}(z)$, as desired.
\end{proof}

Fix $z\in T-T'$ with
a winning strategy $\tau$ for {\it II} in 
${\mathcal G}(z)$.
Define
$\varphi:(\lm)^{<\om}\to\lm$ by
$\varphi(\emptyset)=\mu$
and
$\varphi(s)=e``x\times\mu$, where
$\tau(s)=(b,\pi,x)$.
Since $X_z$ is $F$-positive, 
$\bi_{n<\om}A_n\in X_z$ for some
$\{A_n:n<\om\}\su\lm$ such that
$\varphi(\lan A_k:k<n\ran)\su A_n$ for $n<\om$.
Set
$(b_n,\pi_n,x_n)=\tau(\lan A_k:k\le n\ran)$
for $n<\om$.
Then
$\lan A_n:n<\om\ran$
is a play of {\it I} in ${\mathcal G}(z)$
against $\tau$, since
$\mu=\varphi(\emptyset)\su A_0$
and
$e``x_n\times\mu=\varphi(\lan A_k:k\le n\ran)
\su A_{n+1}$.

Set
$x=\bi_{n<\om}x_n$.
Since
$\{x_n:n<\om\}\su C(f)$ is increasing,
we have
$x\in C(f)$, hence
$e``x\times(x\cap\mu)\su x$.
Also
$x\cap\mu=z\in T$ by $x_n\cap\mu=z$.
Note that
$b_n\in\p_\nu\la$,
$b_n\su x_n=\pi_n``(x_n\cap\mu)\su 
e``b_n\times\mu$ and
$A_n\su e``b_n\times\mu\su 
e``x_n\times\mu\su A_{n+1}$ for $n<\om$.
Hence
$b=\bi_{n<\om}b_n\in\p_\nu x$.
Also
$x\su e``b\times\mu=e``x\times\mu
=\bi_{n<\om}A_n\in X_z=X_{x\cap\mu}$.
Thus
we have $x\in S\cap C(f)$, as desired.
\end{proof}

\begin{clm}
$S$ is $\mu$-nonreflecting.
\end{clm}

\begin{proof}
Suppose to the contrary
$S\cap \p_\ka A$ is
stationary in $\p_\ka A$ for some
$\mu\su A\su\la$ of size $\mu$.
Then 
$\{x\in\p_\ka A:e``x\times(x\cap\mu)\su x\}$ 
is unbounded in $\p_\ka A$,
hence
$e``A\times\mu\su A$.
Moreover
$A=e``a\times\mu$
for some $a\in\p_\nu A$:

Fix a bijection $\pi:\mu\to A$.
Then
$U=\{x\cap\mu:
\pi``(x\cap\mu)=x\in S\cap \p_\ka A\}$
is a stationary subset of $T$.
For $z\in U$ we have $b\in\p_\nu z$ and
$\xi\in z$ such that
$\pi``z\su e``(\pi``z)\times\mu
=e``(\pi``b)\times\mu
\su e``(\pi``c_\xi)\times\mu$, since 
$\pi``z\in S$ and
$\{c_\xi:\xi\in z\}$ is unbounded in $\p_\nu z$. 
Take $\xi<\mu$ and
$U^*\su U$ stationary in $\p_\ka\mu$
so that
$\pi``z\su e``(\pi``c_\xi)\times\mu$ for $z\in U^*$.
Since $\{\pi``z:z\in U^*\}$ 
is stationary in $\p_\ka A$,
$A=\bi_{z\in U^*}\pi``z
\su e``(\pi``c_\xi)\times\mu\su 
e``A\times\mu\su A$.
Hence 
$A=e``(\pi``c_\xi)\times\mu$
and $\pi``c_\xi\in\p_\nu A$, as desired.

For $i=0,1$ take
$a\su x^i\in S\cap\p_\ka A$
so that
$x^i\cap\mu$ disagrees with each other.
Then
$A=e``a\times\mu\su e``x^i\times\mu\su 
e``A\times\mu\su A$.
Hence
$A=e``x^i\times\mu\in X_{x^i\cap\mu}$
by $x^i\in S$,
which contradicts that
$X_{x^i\cap\mu}$ is disjoint
from each other, as desired.
\end{proof}
Therefore $S$ is the desired set.
\end{proof}

Let us derive another

\begin{cor}
$\pkl$ has
a $\ka^+$-nonreflecting stationary subset
if
$\la\ge\ka^{++}$ and
$\cf(\nu,\ga)<\ka$ for $\nu<\ga<\ka$.
\end{cor}

\begin{proof}
It suffices to prove the case $\la=\ka^{++}$
by checking the
conditions of Theorem~3 for
$\la=\mu^+=\ka^{++}$.

For $\ga<\mu=\ka^+$ fix a club set
$T_\ga\su\p_\ka\ga$ of size $\ka$
and for $z\in\bi_{\ga<\mu}T_\ga$
an unbounded set $C_z\su\p_\nu z$
of size $<\ka$.
Set
$\{c_\xi:\xi<\mu\}=
\bi\{C_z:z\in\bi_{\ga<\mu}T_\ga\}$.
Then
$T=\{z\in\bi_{\ga<\mu}T_\ga:
\{c_\xi:\xi\in z\}$
is unbounded in 
$\p_\nu z\}$ has size $\mu$.
We claim that
$T$ is stationary in $\p_\ka\mu$.

Fix $f:[\mu]^{<\om}\to\p_\ka\mu$.
We have 
$\ga<\mu$ of cofinality $\ka$
such that 
$\bi f``[\ga]^{<\om}\cup
\bi_{\xi<\ga}c_\xi\su\ga$ and 
$C_y\su\{c_\xi:\xi<\ga\}$
for
$y\in\bi_{\be<\ga}T_\be$.
Build an increasing and continuous sequence
$\{z_\al:\al<\nu\}\su T_\ga$ so that
$\bi f``[z_\al]^{<\om}\cup
\bi\{c_\xi:\xi\in z_\al\}\su z_{\al+1}$
and $C_y\su\{c_\xi:\xi\in z_{\al+1}\}$
for some
$z_\al\su y\in\bi_{\be<\ga}T_\be$.
Then $z=\bi_{\al<\nu}z_\al\in C(f)$, since
$\bi f``[z_\al]^{<\om}\su z_{\al+1}$.
Since
$\{z_\al:\al<\nu\}\su T_\ga$ is increasing, 
$z\in T_\ga$. 
Since $\bi\{c_\xi:\xi\in z_\al\}\su z_{\al+1}$, 
$\{c_\xi:\xi\in z\}\su\p_\nu z$.
To see that
$\{c_\xi:\xi\in z\}$
is unbounded in $\p_\nu z$, fix $x\in\p_\nu z$.
We have $\al<\nu$ with $x\su z_\al$, hence
$\xi\in z_{\al+1}$ with $x\su c_\xi$, 
as desired.
\end{proof}

Theorem~3 is void, however, if
$\cf\mu<\ka$ or if
$\ka=\theta^+$ and 
$\theta>\cf\theta=\om$:
In the former case
$\p_\ka\mu$
has no stationary subset of size $\mu$.
In the latter case
$\p_\nu z$ has no unbounded subset
of size $\theta$
for $z\in[\mu]^\theta$, since 
$\cf(\nu,\theta)>\theta$ if
$\cf\theta<\nu<\theta$.
See \cite{m} for a nonreflection
result in the latter case 
under additional assumptions.

\section{Proof of Theorem 1}

This section is devoted to 
Foreman--Magidor's example
of a nonreflecting stationary set
as we understand it.
The proof invokes
those~\cite{b, be}
that $\p_\ka\ka^+$
has a club subset of size $\le(\ka^+)^{\om_1}$
and 
that stationary reflection implies
Chang's conjecture.

\begin{proof}[Proof of Theorem 1]
Fix a bijection 
$\pi_\ga:\ka\to\ga$ for
$\ka\le\ga<\ka^+$.
Define 
$h:[\ka^+]^2\to\p_\ka\ka^+$ 
by
$h(\al,\be)=\lim\pi_\be``{\pi_\be}^{-1}(\al)$.
Since $\la\ge2^{\ka^+}$, we have a list
$\{g_\xi:\xi<\la\}$ of the functions
$g:\ka^+\to\p_\ka\ka$.
Then $D=\{x\in\pkl:\bi h``[x\cap\ka^+]^2\su x\land
\fo\ga\in x\cap(\ka^+-\ka) 
(\pi_\ga``(x\cap\ka)=x\cap\ga)\land
\fo\xi\in x(\bi g_\xi``(x\cap\ka^+)\su x)\}$
is club in $\pkl$.

Set
$S=\{x\in\pkl:\{\sup(y\cap\ka^+):
x\su y\in D\land
y\cap\ka=x\cap\ka\}$
is nonstationary in $\ka^+\}$.

\begin{clm}
$S$ is stationary in $\pkl$.
\end{clm}

\begin{proof}
Suppose otherwise.
By induction on $n<\om$ build 
$f_n:[\la]^{<\om}\to\p_\ka\la$
and
$\xi_n:[\la]^{<\om}\to\la$
so that
$C(f_0)\su D-S$,
$g_{\xi_n(a)}(\ga)=\cl_{f_n}(a\cup\{\ga\})\cap\ka$
and 
$f_{n+1}(a)=f_n(a)\cup\{\xi_n(a)\}$.
Define 
$f:[\la]^{<\om}\to\p_\ka\la$ by
$f(a)=\bigcup_{n<\om}f_n(a)$.

\begin{scl}
$\{\sup(z\cap\ka^+):
x\su z\in C(f)\land
z\cap\ka=x\cap\ka\}$
is unbounded in $\ka^+$
for $x\in C(f)$.
\end{scl}

\begin{proof}
Fix $\al<\ka^+$.
Since $x\in C(f)\su\pkl-S$, 
$\{\sup(y\cap\ka^+):
x\su y\in D\land
y\cap\ka=x\cap\ka\}$
is stationary in $\ka^+$.
Hence we have
$x\su y\in D$ with
$y\cap\ka=x\cap\ka$ and
$\al<\ga\in y\cap\ka^+$.

Set
$z=\bi\{\cl_{f_n}(a\cup\{\ga\}):
n<\om\land a\in[x]^{<\om}\}$.
Then $\al<\ga\le\sup(z\cap\ka^+)$.
It is easy to see
$x\su z\in C(f)$.
To see
$z\cap\ka\su x\cap\ka$,
fix
$\be\in z\cap\ka$.
Then
$\be\in\cl_{f_n}(a\cup\{\ga\})\cap\ka=
g_{\xi_n(a)}(\ga)$
for some $n<\om$ and $a\in[x]^{<\om}$.
Since $x\in C(f)$ and $a\in[x]^{<\om}$,
$\xi_n(a)\in f(a)\su x\su y$.
Hence
$\be\in g_{\xi_n(a)}(\ga)\su
y\cap\ka=x\cap\ka$, 
as desired, since
$\xi_n(a),\ga\in y\in D$.
\end{proof}

For $i=0,1$
build an increasing and continuous sequence
$\{x^i_\xi:\xi<\om_1\}\su C(f)$
so that
$x^i_\xi\cap\ka=x^0_0\cap\ka\in\ka$
has cofinality $\om_1$,
$\sup(x^0_\xi\cap\ka^+)\le
\sup(x^1_\xi\cap\ka^+)
<\sup(x^0_{\xi+1}\cap\ka^+)$ and
$x^1_0\cap\ka^+$ is not an initial segment of
$x^0_1\cap\ka^+$
as follows:
First take $x^0_0\in C(f)$ with
$x^0_0\cap\ka\in S^{\om_1}_\ka$.
Subclaim allows us to take $x^0_1$ from
$X=\{z\in C(f):
x^0_0\su z\land
z\cap\ka=x^0_0\cap\ka\}$
so that
$\{\sup(z\cap\ka^+):z\in X\}\cap
\sup(x^0_1\cap\ka^+)$
has size $\ka$.
Since $x^0_1\cap\ka^+$
has $<\ka$ initial segments,
we have $x^1_0\in X$ as required above.
The rest of the construction is routine.

Set $x^i=\bi_{\xi<\om_1}x^i_\xi$.
Then
$x^i\in C(f)$, since
$\ka\ge\om_2$ is regular and
$\{x^i_\xi:\xi<\om_1\}\su C(f)$
is increasing.
Also
$\sup(x^i\cap\ka^+)=
\sup_{\xi<\om_1}\sup(x^i_\xi\cap\ka^+)$
has cofinality $\om_1$
and agrees with each other by
$\sup(x^0_\xi\cap\ka^+)\le
\sup(x^1_\xi\cap\ka^+)
<\sup(x^0_{\xi+1}\cap\ka^+)$.
Since $x^i,x^i_\xi\in C(f)\su D$,
we have
$x^i\cap\ga=
\pi_\ga``(x^i\cap\ka)=\pi_\ga``(x^0_0\cap\ka)=
\pi_\ga``(x^i_\xi\cap\ka)=x^i_\xi\cap\ga$
for
$\ga\in x^i_\xi\cap(\ka^+-\ka)$.
Since
$x^1_0\cap\ka^+$ is not an initial segment of
$x^0_1\cap\ka^+$,
$x^i\cap\ka^+$ disagrees
with each other.
Moreover
$x^i\cap\ka^+$ is $\si$-closed:

Fix $b\su x^i\cap\ka^+$
of order type $\om$.
We have $b\su\be\in x^i\cap(\ka^+-\ka)$
by
$\cf\sup(x^i\cap\ka^+)=\om_1$.
Since 
${\pi_\be}^{-1}``(x^i\cap\be)
=x^i\cap\ka=x^0_0\cap\ka\in\ka$
has cofinality $\om_1$,
we have
$\al\in x^i\cap\be$ with
${\pi_\be}^{-1}``b\su{\pi_\be}^{-1}(\al)$.
Hence 
$b\su\pi_\be``{\pi_\be}^{-1}(\al)$.
Thus
$\sup b\in h(\al,\be)\su x^i$, as desired,
since $\al,\be\in x^i\in D$.

Set
$c=x^0\cap x^1\cap\ka^+$, which
is unbounded in
$\sup(x^i\cap\ka^+)$.
Then
$x^i\cap\ka^+=\bi_{\ga\in c}x^i\cap\ga
=\bi_{\ga\in c}\pi_\ga``(x^i\cap\ka)
=\bi_{\ga\in c}\pi_\ga``(x^0_0\cap\ka)$
by $x^i\in D$,
which contradicts that
$x^i\cap\ka^+$
disagrees with each other, as desired.
\end{proof}

\begin{clm}
$S$ is nonreflecting.
\end{clm}

\begin{proof}
Suppose to the contrary
$S\cap\p_\ka A$ is stationary
in $\p_\ka A$
for some $\ka\su A\su\la$
of size $\ka$.
Fix a bijection $\pi:\ka\to A$.
Then
$T=\{\ga<\ka:\pi``\ga\in S\land
\pi``\ga\cap\ka=\ga\}$ 
is stationary in $\ka$,
hence
$\{y\cap\ka^+:
\pi``(y\cap\ka)\su y\in D
\land y\cap\ka\in T\}$
is stationary in $\p_\ka\ka^+$.
Thus
$\{\sup(y\cap\ka^+):
\pi``(y\cap\ka)\su y\in D
\land y\cap\ka\in T\}$
is stationary in $\ka^+$,
hence so is
$\{\sup(y\cap\ka^+):
\pi``(y\cap\ka)\su y\in D
\land y\cap\ka=\ga\}$
for some $\ga\in T$.
Thus
$\{\sup(y\cap\ka^+):
\pi``\ga\su y\in D
\land y\cap\ka=\pi``\ga\cap\ka\}$
is stationary in $\ka^+$,
which contradicts
$\pi``\ga\in S$, as desired.
\end{proof}
Therefore stationary reflection for $\pkl$ fails.
\end{proof}

We remark that the same proof as above 
works if
we replace ``nonstationary" by
``bounded"
in the above definition of $S$.

\end{document}